\documentclass{amsproc}
\usepackage{graphicx} 
\usepackage{xcolor, tikz}
\usepackage{amsmath,latexsym,amssymb,amsthm,enumerate,amscd}
\usepackage{arydshln}

\newtheorem{theorem}{Theorem}[section]
\newtheorem{lemma}[theorem]{Lemma}
\newtheorem{proposition}[theorem]{Proposition}

\newtheorem{definition}[theorem]{Definition}
\newtheorem{example}[theorem]{Example}
\newtheorem{conjecture}[theorem]{Conjecture}

\numberwithin{equation}{section}

\newcommand\scalemath[2]{\scalebox{#1}{\mbox{\ensuremath{\displaystyle #2}}}}

\def\k{\mathsf{k}}

\def\rank{\mathrm{rank}}

\def\<{\left<}
\def\>{\right>}

\begin{document}
\title{Commuting Jordan Types: a Survey}
\author{Leila Khatami} 
\address{Department of mathematics, Union College, Schenectady, NY 12308, USA}
\curraddr{}
\email{khatamil@union.edu}
\thanks{}


\date{}

\begin{abstract}
In this paper, we survey the progress in the problem of finding the maximum commuting nilpotent orbit that intersects the centralizer of a given nilpotent matrix.
 
\end{abstract}
\maketitle

\section{Introduction}
Let $\k$ be an infinite field and $n$ be a positive integer. We denote the set of all $n \times n$ matrices with entries in $\k$ by ${\mathcal Mat}_n(\k)$, and the set of nilpotent matrices in ${\mathcal Mat}_n(\k)$ by ${\mathcal Nilp}_n(\k)$. For a matrix  $A \in {\mathcal Nilp}_n(\k)$, let $\mathcal{N}_A$ be the nilpotent commutator of $A$, that is
$
\mathcal{N}_A=\{B \in {\mathcal Nilp}_n(\k) \, | \, [A,B]=0 \}.
$

By Jordan decomposition theorem, there is a bijection between the set of partitions of $n$ and $\mbox{GL}_n(\k)$-orbits of ${\mathcal Nilp}_n(\k)$. For a nilpotent matrix  $A\in {\mathcal Nilp}_n(\k)$, we refer to the partition of $n$ determined by the size of Jordan blocks in the Jordan canonical form of $A$ as the {\it Jordan type} of $A$ and we denote it by $P_A$. It is well known that the nilpotent commutator of a nilpotent matrix is an irreducible algebraic variety (for example see \cite[Proposition 2.3]{Basili03}). This in particular implies that the nilpotent commutator $\mathcal{N}_A$ contains an open dense subset over which the Jordan type is constant. In the past two decades, many authors have studied the characteristics and properties of the generic Jordan type in the nilpotent commutator of a nilpotent matrix. In this paper, we  provide a general overview of those results. A wide variety of algebraic and combinatorial tools are used to establish the results that we survey in this paper. In the interest of brevity, we leave out technical details and keep the focus of this paper on the results themselves.  

We should note that the study of commuting nilpotent matrices in the past two decades has not been limited to the results surveyed in this paper. For example, see the work of Baranovsky, Basili and Iarrobino in \cite{Bar, Basili03, BI} where the study of pairs of commuting nilpotent matrices has been related to the punctual Hilbert scheme of a plane. Pairs of commuting nilpotent orbits have also been studied by Premet, McNinch and Panyushev in the broader context of Lie algebras (see \cite{McNinch, Pa, Premet}). Triples (or more generally $d$-tuples) of commuting nilpotent matrices have also been studied by Guralnick, Sethuranam, \v{S}ivic, Haboush and Hyeon, and others (see \cite{GS, HaboushHyeon, NeSe, NS, SeSi, S, Sivic3, Sivic2}).

\section{Preliminaries}
Throughout this paper, we write a partition as a weakly decreasing sequence of positive integers. If a partition includes $k$ parts of size $\ell$, in writing $P$, we often include $k^\ell$ (instead of writing $\ell$ copies of $k$). For example, we write $(3^2, 1^3)$ for $(3,3,1,1,1)$. 
\begin{definition}\label{commdefn}

    Let $n$ be a positive integer and $P$ and $Q$ be partitions of $n$. We say partitions $P$ and $Q$ {\it commute} if there exist commuting nilpotent matrices $A$ and $B$ in ${\mathcal Nilp}_n(\k)$ with Jordan types $P$ and $Q$, respectively.  
\end{definition}

The Jordan type of a nilpotent matrix can be determined by the ranks of its powers.
\begin{lemma}(\cite[Lemma~1.1]{BI})
Suppose that $A\in {\mathcal Nilp}_n(\k)$. For $i\geq 1$, let $d_i=\operatorname{rank} A^{i-1}-\operatorname{rank} A^{i}$. Then the partition $d=(d_1, d_2, \dots)$ and the Jordan type $P_A$ of $A$ are conjugate partitions.
\end{lemma}
\begin{example}\label{4_22ex}
    Let $J_4$ be the $4 \times 4$ nilpotent Jordan block. Obviously the Jordan type of $J_P$ is $(4)$. Now consider $B=J_4^2$. It is easy to see that $\operatorname{rank} B=2$, while $B^2=0$. Thus, the Jordan type of $B$ is a conjugate of $(2,2)$, which is incidentally $(2,2)$ itself.  This in particular shows that partitions $(4)$ and $(2,2)$ commute. Note that although partitions $(2,2)$ and $(4)$ commute, Jordan matrix $J_{(2,2)}$, the nilpotent Jordan matrix corresponding to the partition $(2,2)$, and $J_4$ do not commute.
\end{example}

Suppose that matrices $A$ and $B$ commute. Then for every matrix $A'$ in the $\operatorname{GL}_n(\k)$-orbit of $A$, there exists a matrix $B'$ in the $\operatorname{GL}_n(\k)$-orbit of $B$ such that $A'$ and $B'$ commute (if $A'=P^{-1}AP$, simply let $B'=P^{-1}BP$). However, as the example above also shows, this does not imply that every matrix in the orbit of $A$ must commute with all matrices in the orbit of $B$.

We can extend the approach we took in Example~\ref{4_22ex} to easily find Jordan types of powers of a single Jordan block. This will immediately give rise to a family of partitions of $n$ that commute with partition $(n)$. 

Let $J_n$ be the $n \times n$ nilpotent Jordan block. Then multiplying any $n \times n$ matrix $A$ by $J_n$ shifts the components of each column of $A$ one position up, with a zero in the last position  (we work with upper triangular Jordan blocks). In particular, for $1\leq k \leq n$, the matrix $J_n^k$ has 1's on its $k$-diagonal and 0's everywhere else. Consequently, 

$$\operatorname{rank} J_n^k=\left\{\begin{array}{lll}n-k,&&1\leq k\leq n,\\ 0,&&\mbox{otherwise.}\end{array}\right.$$ Therefore, if non-negative integers $x$ and $r$ are such that $0\leq r<k$ and $n=kx+r$, then $(\operatorname{rank}\left(J_n^k\right)^i)_{i\geq 0}=(n,n-k,n-2k,\dots,n-kx,0)$. Thus the difference sequence, which is conjugate to the Jordan type of $J_n^k$ is $(k^x,r)$. Therefore, the Jordan type of $J_n^k$ is the partition $\Big((x+1)^r,x\Big)$. Note that this is the unique partition of $n$ consisting of exactly $k$ ``almost equal" parts.
\begin{definition}\label{arpartsdefn}
    A partition is called {\it almost rectangular} if its largest and smallest parts differ by at most 1. For a positive integer $n$ and $1\leq k \leq n$, the almost rectangular partition of $n$ with $k$ parts is denoted by $[n]^k$. If $r$ denotes the remainder of dividing $n$ by $k$, then the almost rectangular partition $[n]^k$ consists of $r$ (possibly 0) parts of size $\lfloor\frac{n}{k}\rfloor+1$ and $k-r$ parts of size $\lfloor\frac{n}{k}\rfloor$.

    Every partition $P$ can be decomposed into a union (concatenation) of its almost rectangular subpartitions. We denote the minimum number of non overlapping almost rectangular subpartitions required for such a decomposition of $P$ by $r_P$.
\end{definition}
\begin{example}
    Partitions $[6]^2=(3,3)$ and $[6]^4=(2^2,1^2)$ are both almost rectangular partitions of $n=6$. 
    
    Now consider $P=(5,3,2^2,1^2)$. Then we can write 
    $$P=(\underline{5},\underline{3,2^2},\underline{1^2})=([5]^1,[7]^3,[2]^2) \mbox{ or }P=(\underline{5},\underline{3}, \underline{2^2,1^2})=([5]^1,[3]^1,[6]^4).$$ 
    The partition can not be decomposed into fewer almost rectangular subpartitions. Therefore, $r_P=3$.
\end{example}

\section{The generic commuting partition}
\begin{definition}\label{mapQdefn}
     Let $n$ be a positive integer and $\mathcal{P}(n)$ be the set of all partitions of $n$. We define the map $\mathcal{Q}:\mathcal{P}(n) \to \mathcal{P}(n)$ to send each partition $P$ to the Jordan type of a generic element in the irreducible variety $\mathcal{N}_{J_P}$, where $J_P$ is the nilpotent Jordan matrix with Jordan type $P$. In other words, $\mathcal{Q}(P)$ is the generic Jordan type that commutes with $P$.
\end{definition}

In this section, we give an overview the properties of the map $\mathcal{Q}$ above.

Recall that the {\it dominance order} is a partial order defined on the set of all partitions of $n$ as follows.  For partitions $P=(p_1, p_2, \dots)$ and $Q=(q_1, q_2, \dots)$ in $\mathcal{P}(n)$, we say $Q\geq P$ if and only if for all $k$, $$\sum\limits_{1\leq i \leq k}q_i\geq \sum\limits_{1\leq i \leq k}p_i.$$ Equivalently, $Q\geq P$ if and only if for all $k$, $\rank\, {J_Q}^k\geq \rank \, {J_P}^k$. 

Using semicontinuity of the ranks of powers of a matrix and the definition of $\mathcal{Q}(P)$ as the Jordan type of a generic element in the nilpotent commutator of $J_P$, we get the following statement.
\begin{proposition}
    Let $n$ be a positive integer and $P$ be a partition of $n$. If partition $Q$ commutes with $P$, then $\mathcal{Q}(P)\geq Q$. 
\end{proposition}
In other words, $Q(P)$ is the maximum partition (in the dominance order) among Jordan types of matrices in the nilpotent commutator of the Jordan matrix $J_P$ (or any other matrix with Jordan type $P$). 

In \cite{BI}, Basili and Iarrobino characterize partitions that are stable under the action of $\mathcal{Q}$. In \cite[Theorem 2.1]{Pa}, Panyushev independently proved the same result in a more general context. 
\begin{theorem}(\cite[Theorem~1.12]{BI})\label{stable}
    For a partition $P$, $\mathcal{Q}(P)=P$ if and only if parts of $P$ differ pairwise by at least 2. 
\end{theorem}

 For a pair of commuting nilpotent matrices $A$ and $B$ in ${\mathcal Nilp}_n(\k)$, consider the subalgebra $\k[A,B]$ of ${\mathcal Mat}_n(\k)$ that is generated by $A$ and $B$. Then $\k[A,B]$ is a commutative Artinian $k$-algebra. Let $H(\k[A,B])$ denote the Hilbert function of this algebra and assume that $\operatorname{char} \k$ is either $0$ or is greater than the socle degree of $H(\k[A,B])$. In \cite[Theorem 2.16]{BI}, Basili and Iarrobino prove that if $\dim_{\k}\k[A,B]=n$, then the Jordan type of the multiplication by the general linear form in $\k[A,B]$, as well as its associated graded algebra, is the conjugate of the partition obtained from $H(\k[A,B])$. On the other hand, in \cite[Corollary 5]{KO}, Ko\v{s}ir and Oblak show that if $B$ is a generic element in $\mathcal{N}_A$, then the algebra $\k[A,B]$ is Gorenstein and therefore also complete intersection. Using Macaulay's classical characterization of Hilbert functions of complete intersection Artinian algebras, Ko\v{s}ir and Oblak then prove the following statement.
\begin{theorem}(\cite[Theorem 6]{KO})\label{idempotent}
    Assume that $\operatorname{char} \k=0$ or $\operatorname{char} \k >n$. Then the map $\mathcal{Q}$ defined in Definition~\ref{mapQdefn} is idempotent, i.e., for every partition $P$, $\mathcal{Q}(\mathcal{Q}(P))=P$.
\end{theorem}

In \cite{Basili03}, Basili determines the number of parts of $\mathcal{Q}(P)$ by calculating the rank of a generic matrix in the nilpotent commutator of $J_P$. Also see \cite[Theorem~2.17]{BIK} for another proof. 

\begin{theorem}(\cite[Proposition 2.4]{Basili03})
    Partition $\mathcal{Q}(P)$ has $r_P$ parts. 
\end{theorem}
Note that an immediate consequence of this result is the fact that a partition $P$ of $n$ commutes with the maximal partition $(n)$ if and only if $P$ is almost rectangular. 

In \cite[Theorem 12]{Oblak}, Oblak gives an explicit formula for the largest part of $\mathcal{Q}(P)$. Although her result was proved over a field of characteristic 0, it was later generalized to an arbitrary infinite field by Iarrobino and the author in \cite[Corollary 3.10]{IK}. In \cite{K:smallest}, the author gives an explicit formula for the smallest part of $\mathcal{Q}(P)$ over an arbitrary infinite field. 

These partial results in particular determined the partition $\mathcal{Q}(P)$ for every partition $P$ with $r_P\leq 3$. They also confirmed a conjecture of Oblak that was finally proved by Basili in \cite{basili21} in its generality after being partially proved in \cite{IK} and \cite{K:unique}. These results are closely tied to a recursive process introduced by Oblak. In the next section, we discuss the process, as well as an algorithm to explicitly determine $\mathcal{Q}(P)$ from $P$.

\section{Oblak Process}

In this section, we review a recursive process, proposed by Oblak, which she conjectured could be used to explicitly determine the generic commuting partition $\mathcal{Q}(P)$ in terms of $P$. As we noted earlier, this conjecture was partially confirmed through the work of several authors and was finally proved in general by Basili in \cite{basili21}. The material presented in this section is based on Oblak's original ideas but the notations and the set-up are mainly based on \cite{BIK}, \cite{IK}, \cite{K:smallest} and \cite{K:unique}. 

\noindent{\bf The poset $\mathcal{D}_P$.} Let $n$ be a positive integer and $P$ be a partition of $n$. We write $P=(p_1^{n_1}, p_2^{n_2}, \dots, p_t^{n_t})$, where $p_1>\dots>p_t$ and $n_i>0$, for all $i$. We define the directed graph $\mathcal{D}_P$ as follows. 

There are $n$ vertices in $\mathcal{D}_P$. The vertices are arranged in rows such that each row corresponds to a part of $P$. The rows are arranged in a weakly decreasing order from top to bottom of $\mathcal{D}_P$. In other words, there are $n_1$ rows of length $p_1$ followed by $n_2$ rows of length $p_2$ below them, etc. We label each vertex by a triple $(u,p_i,k)$, where for each $i$, $1\leq u\leq p_i$ and $1\leq k \leq n_i$. The vertex $(u,p_i,k)$ corresponds to the $u$-th vertex in the $k$-th length-$p_i$ row of $\mathcal{D}_P$. For each $i$, the vertices in the $n_i$ rows of length $p_i$ are labeled in a way that the first (respectively, last) component of the triple is increasing when we go from left to right (from top to bottom, respectively). 

\begin{figure}[ht]

\begin{tikzpicture}[scale=0.9, transform shape]
\newcommand\Square[1]{+(-#1,-#1) rectangle +(#1,#1)}
\foreach \x in {1,2,3,4,5,6}{\filldraw (\x,10-1) circle (2.5pt);}
\foreach \x in {1,2,3,4}{\filldraw (\x+1, 10-3) circle (2.5pt);}
\foreach \x in {1,2,3,4}{\filldraw (\x+1,10-4) circle (2.5pt);}
\foreach \x in {1,2,3}{\filldraw (\x+1.5,10-5) circle (2.5pt);}
\foreach \x in {1,2,3}{\filldraw (\x+1.5,10-6) circle (2.5pt);}
\foreach \x in {1,2}{\filldraw (\x+2,10-7) circle (2.5pt);}
\foreach \x in {1,2}{\filldraw (\x+2,10-8) circle (2.5pt);}
\foreach \x in {1}{\filldraw (\x+2.5,10-9) circle (2.5pt);}

\foreach \x in {1,2,3,4,5,6}{\draw node at(\x,10-1+0.25) {\tiny{$(\x,6,1)$}};}
\draw node at(1+2-0.55,10-7) {\tiny{$(1,2,1)$}};
\draw node at(2+2+0.55,10-7) {\tiny{$(2,2,1)$}};
\draw node at(1+2-0.5,10-8) {\tiny{$(1,2,2)$}};
\draw node at(2+2+0.5,10-8) {\tiny{$(2,2,2)$}};

\foreach \x in {1,2,3,4,5} {\draw [dashdotted, -stealth](\x,10-1)--(\x+0.9,10-1);}
\foreach \x in {1,2,3,4} {\draw [dashed, -stealth](\x,10-1)--(\x+0.95,10-3.9);}
\foreach \x in {1,2,3} {\draw [dashed, -stealth](\x+1+0.05,10-3)--(\x+1.5-0.05,10-5.9);}
\foreach \x in {1,2} {\draw [dashed,-stealth](\x+1.5+0.05,10-5)--(\x+2-0.05,10-7.9);}
\foreach \x in {1} {\draw [dashed,-stealth](\x+2+0.05,10-7)--(\x+2.5-0.05,10-8.9);}
\foreach \x in {1} {\draw [densely dotted, stealth-](\x+2+1+0.05,10-8.1)--(\x+2.5,10-9);}
\foreach \x in {1,2} {\draw [densely dotted, stealth-](\x+1.5+1,10-6.1)--(\x+2,10-7);}
\foreach \x in {1,2,3} {\draw [densely dotted, stealth-](\x+1+1,10-4.1)--(\x+1.5,10-5);}
\foreach \x in {1,2,3,4} {\draw [densely dotted, stealth-](\x+2,10-1.1)--(\x+1,10-3);}
\foreach \x in {1,2,3,4} {\draw [-stealth](\x+1,10-4)--(\x+1,10-3.1);}
\foreach \x in {1,2,3} {\draw [-stealth](\x+1.5,10-6)--(\x+1.5,10-5.1);}
\foreach \x in {1,2} {\draw [-stealth](\x+2,10-8)--(\x+2,10-7.1);}
\end{tikzpicture}

\caption{The directed graph $\mathcal{D}_P$ corresponding to the partition $P=(6,4^2,3^2,2^2,1)$.}\label{poset}
\end{figure}
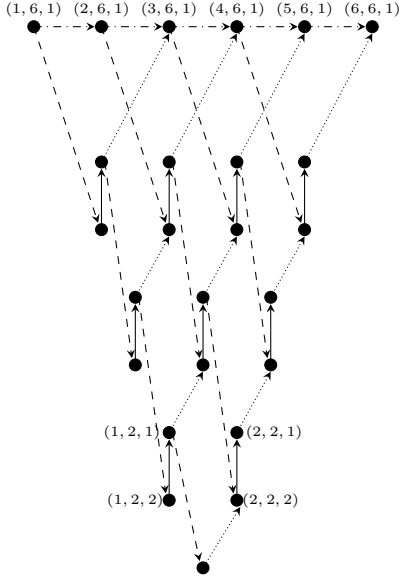

For $1\leq i< t$, there are four families of edges in $\mathcal{D}_P$. 
\begin{itemize}
    \item For each $1\leq u \leq p_{i+1}$, there is an edge from $(u,p_i, 1)$, the $u$-th vertex in the top row corresponding to $p_i$, to $(u, p_{i+1}, n_i)$, the $u$-th vertex in the bottom row corresponding to $p_{i+1}$ (dashed arrows in Figure~\ref{poset}).
    \item For each $1\leq u \leq p_{i+1}$, there is edge from $(u,p_{i+1},1)$, the $u$-th vertex in the top row corresponding to $p_{i+1}$ to $(u+p_i-p_{i+1},p_i,n_i)$, the $u+(p_i-p_{i+1})$-th vertex in the bottom row corresponding to $p_{i}$ (dotted arrows in Figure~\ref{poset}).
    \item For each $1\leq u \leq p_{i}$ and each $1\leq k <n_i$, there is an edge from $(u,p_i,k)$, the $u$-th vertex in the $k$-th row corresponding to $p_i$, to $(u,p_i,k-1)$ the $u$-th vertex in the $(k-1)$-st row corresponding to $p_{i}$ (solid arrows in Figure~\ref{poset}).   
    \item If $p_i$ is an isolated part, i.e, $p_{i-1}-p_i>1$ and $p_{i}-p_{i+1}>1$, then for every $1\leq u<p_i$, the directed graph $\mathcal{D}_P$ also includes an edge from $(u,p_i,1)$, the $u$-th vertex in the top row corresponding to $p_i$, to $(u+1,p_i,n_i)$, the $(u+1)$-st vertex in the bottom row corresponding to $p_i$ (dash-dotted arrows in arrows in Figure~\ref{poset}). Note that if $n_i=1$, then these edges go from left to right in the same row (the only row of $\mathcal{D}_P$ corresponding to $p_i$).
\end{itemize}

In Figure~\ref{poset}, we illustrate the poset $\mathcal{D}_{P}$ for $P=(6,4^2,3^2,2^2,1)$. We have only labeled a few vertices of the poset by their corresponding triples to illustrate the labeling logic without avoid overcrowding the figure.

The directed graph $\mathcal{D}_P$ is in fact the covering edge graph of a partially ordered set (poset), which we also denote by $\mathcal{D}_P$. For elements $v$ and $v'$ in the poset, $v>v'$ if and only if there is a directed path from vertex $v$ to vertex $v'$ in the directed graph $\mathcal{D}_P$. The partial order is defined over a basis of $\k^n$ in close connection with a maximal nilpotent subalgebra of the centralizer of $J_P$. The generic Jordan type in this algebra is $\mathcal{Q}(P)$, as well. We refer the reader to \cite{BIK} and \cite{K:unique} for details. Roughly speaking, the edges in $\mathcal{D}_P$ correspond to non-zero entries in a generic (upper triangular) matrix in the nilpotent commutator of $J_P$.  As a simple example, consider partition $P=(3,1)$. Straightforward calculations show that a generic element in the nilpotent commutator of $J_P$ is a matrix $B$ of the form given in Figure~\ref{(3,1)}, where $a,b,c$ are general enough. The Figure also illustrates $\mathcal{D}_P$ with each edge labeled by the corresponding entry in $B$.
\begin{figure}
    \centering
\begin{tikzpicture}[scale=1, transform shape]
\draw[black]
node at (0,-2){$B=\scalemath{1}{\begin{pmatrix}0&a&0&|b\\0&0&a&|0\\0&0&0&|0\\ \hdashline 0&0&c&|0\end{pmatrix}}$ }
node at (0,8-12){Generic matrix in the }
node at (0,8-12.5){Nilpoetent commutator of $J_{(3,1)}$};

\foreach \x in {1,2,3}{\filldraw (\x+5,10-1-10) circle (2.5pt);}
\foreach \x in {1}{\filldraw (\x+1+5,10-3-10) circle (2.5pt);}

\foreach \x in {1,2} {\draw [dashdotted, -stealth](\x+5,10-1-10)--(\x+0.9+5,10-1-10);}
\foreach \x in {1} {\draw [dashed, -stealth](\x+5,10-1-10)--(\x+0.95+5,10-2.9-10);}
\foreach \x in {1} {\draw [densely dotted, stealth-](\x+2+5,10-1.1-10)--(\x+0.95+5,10-3-10);}

\draw[black]
node at (1.5+5,9.25-10){$a$}
node at (2.5+5,9.25-10){$a$}
node at (1.25+5,8-10){$b$}
node at (2.75+5,8-10){$c$}
node at (2+5,8-12){$\mathcal{D}_{(3,1)}$};
\end{tikzpicture}
    
    \caption{A generic matrix in the nilpotent commutator of $J_P$, as well as the covering edge graph of poset $\mathcal{D}_P$ for $P=(3,1)$ and the poset }
    \label{(3,1)}
\end{figure}

Recall that a totally ordered subset of a poset is called a {\it chain}. In \cite{Gansner}, Gansner proves that in a finite poset $S$, the Jordan type of a generic element in the incidence algebra of $S$, is obtained by the difference sequence of the sequence of $\{c_k\}_{k}$, where $c_k$ is the maximal cordiality of a union of $k$ chains in the poset. We denote this partition by $\lambda(S)$ (see \cite{Britz-Fomin}).

As our example above illustrates, each element of the nilpotent commutator of $J_P$ is a matrix in the incidence algebra of $\mathcal{D}_P$. However,  additional conditions on entries of the matrix must be imposed in order for the commuting condition is satisfied. Consequently, the partition $\mathcal{Q}(P)$ is dominated by $\lambda(\mathcal{D}_P)$ (see \cite{IK}). Oblak's recursive process gives an algorithm to determine $\mathcal{Q}(P)$ from the poset $\mathcal{D}_P$ by studying special chains in the poset.  These chains correspond to almost rectangular subpartitions of $P$. 
\begin{definition}
    Let $n$ be a positive integer and $P$ be a partition of $n$. For a positive integer $p$, let $N_P(p)$ denote the (potentially zero) number of parts of $P$ of size $p$. Suppose that positive integer $p$ is such that $N_P(p)>0$, i.e., $P$ has at least one part of size $p$. Let $R_P(p)=\left(p^{N_P(p)}, (p-1)^{N_P(p-1)}\right)$ be the almost rectangular subpartition of $P$ with greatest part of size $p$. Note that $P$ may not include any parts of size $p-1$. If that is the case, then $N_P(p-1)=0$ and $R_P(p)$ will consist of $N_P(p)$ parts of size $p$. 
\end{definition}  

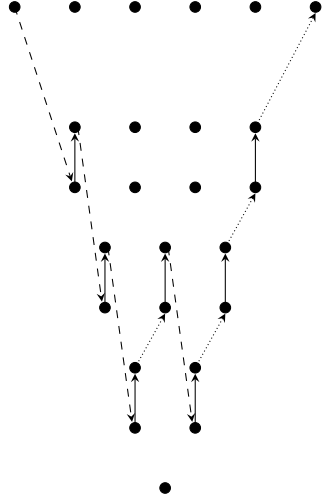
\begin{figure}[ht]

\begin{tikzpicture}[scale=0.8, transform shape]
\newcommand\Square[1]{+(-#1,-#1) rectangle +(#1,#1)}
\foreach \x in {1,2,3,4,5,6}{\filldraw (\x,10-1) circle (2.5pt);}
\foreach \x in {1,2,3,4}{\filldraw (\x+1, 10-3) circle (2.5pt);}
\foreach \x in {1,2,3,4}{\filldraw (\x+1,10-4) circle (2.5pt);}
\foreach \x in {1,2,3}{\filldraw (\x+1.5,10-5) circle (2.5pt);}
\foreach \x in {1,2,3}{\filldraw (\x+1.5,10-6) circle (2.5pt);}
\foreach \x in {1,2}{\filldraw (\x+2,10-7) circle (2.5pt);}
\foreach \x in {1,2}{\filldraw (\x+2,10-8) circle (2.5pt);}
\foreach \x in {1}{\filldraw (\x+2.5,10-9) circle (2.5pt);}

\foreach \x in {1} {\draw [dashed, -stealth](\x,10-1)--(\x+0.95,10-3.9);}
\foreach \x in {1} {\draw [dashed, -stealth](\x+1+0.05,10-3)--(\x+1.5-0.05,10-5.9);}
\foreach \x in {1,2} {\draw [dashed,-stealth](\x+1.5+0.05,10-5)--(\x+2-0.05,10-7.9);}
\foreach \x in {1,2} {\draw [densely dotted, stealth-](\x+1.5+1,10-6.1)--(\x+2,10-7);}
\foreach \x in {3} {\draw [densely dotted, stealth-](\x+1+1,10-4.1)--(\x+1.5,10-5);}
\foreach \x in {4} {\draw [densely dotted, stealth-](\x+2,10-1.1)--(\x+1,10-3);}
\foreach \x in {1,4} {\draw [-stealth](\x+1,10-4)--(\x+1,10-3.1);}
\foreach \x in {1,2,3} {\draw [-stealth](\x+1.5,10-6)--(\x+1.5,10-5.1);}
\foreach \x in {1,2} {\draw [-stealth](\x+2,10-8)--(\x+2,10-7.1);}

\end{tikzpicture}

\caption{The $U$-chain $U_P(3)$ corresponding to the almost rectangular subpartition $R_P(3)=(3^2,2^2)$ of $P=(6,4^2,3^2,2^2,1)$.}\label{Uchainfig}
\end{figure}

In $\mathcal{D}_P$, we define a $U$-chain to be a chain that consists of all vertices in $\mathcal{D}_P$ that correspond to an almost rectangular subpartition of $P$, as well as the first and the last vertex in any row above that almost rectangular subpartition (see Figure~\ref{Uchainfig} for an example). For a positive integer $p$, we refer to the $U$-chain corresponding to the almost rectangular subpartition $R_P(p)$ by $U_P(p)$. In other words, $U_P(p)$ is defined as the set consisting of all vertices of $\mathcal{D}_P$ corresponding to $R_P(p)$, as well as the first and the last vertex in every row of $\mathcal{D}_P$ that is longer than $p$ (see Figure~\ref{Uchainfig} and Example~\ref{Uchainexa}).

Evidently, the cardinality of $U_P(p)$, which we denote by $\mathfrak{u}_P(p)$, is \break
    $\mathfrak{u}_P(p)=pN_P(p)+(p-1)N_P(p-1)+2\sum\limits_{j>p}N_P(j).$ Here the summand $pN_P(p)+(p-1)N_P(p-1)$ is the size of the almost rectangular subpartition $R_P(p)$, while $\sum\limits_{j>p}N_P(j)$ is the number of parts of $P$ that are greater than $p$.
\begin{example}\label{Uchainexa}
    Figure~\ref{Uchainfig}, illustrates the $U$-chain $U_P(3)$ for $P=(6,4^2,3^2,2^2,1)$. Here the almost rectangular subpartition $R_P(3)$ is $(3^2,2^2)$ and $$\mathfrak{u}_P(3)=3(2)+2(2)+2(2+1)=16.$$

    For the same partition $P=(6,4^2,3^2,2^2,1)$, there are other $U$-chains in $\mathcal{D}_P$ corresponding to almost rectangular subpartitions $R_P(6)=(6)$, $R_P(4)=(4^2,3^2)$, $R_P(2)=(2^2,1)$, and even $R_P(1)=(1)$. We have 
    $$\begin{array}{ll}
    \mathfrak{u}_P(6)=&6(1)=6,\\
    \mathfrak{u}_P(4)=&4(2)+3(2)+2(1)=16,\\
    \mathfrak{u}_P(3)=&3(2)+2(2)+2(2+1)=16,\\
    \mathfrak{u}_P(2)=&2(2)+1(1)+2(2+2+1)=15, \mbox{ and}\\
    \mathfrak{u}_P(1)=&1(1)+2(2+2+2+1)=15.
    \end{array}$$
\end{example}

\noindent{\bf Oblak Process.} 
The Oblak process for a partition $P$ is a recursive process to obtain $\mathcal{Q}(P)$ by finding a $U$-chain in $\mathcal{D}_P$ with maximum possible cardinality, removing it from $\mathcal{D}_P$ to obtain a new partition, and repeating the process until all vertices of $\mathcal{D}_P$ are accounted for. 
 
Let $P=(p_1^{n_1}, \dots, p_t^{n_t})$ such that $p_1>\dots>p_t$ and $n_i>0$ for $1\leq i\leq t$.  Assume that $p\in \{p_1, \dots, p_t\}$. Consider the partition $P_1$ determined by the number of remaining vertices in $\mathcal{D}_P-U_P(p)$. In other words, $P_1$ is obtained from $P$ by removing the almost rectangular subpartition $R_P(p)$, keeping the subpartition of $P$ that is below $R_P(p)$ unchanged, and shortening each part of $P$ that is above $R_P(p)$ by 2. If we write $P_1=(q_1^{m_1}, \dots, q_t^{m_t})$, then for $1\leq i \leq t$, we have
 $$\begin{array}{lll}
 
      m_i=\left\{\begin{array}{ll}
     0,&\mbox{if } p_i\in\{p,p-1\},\\
     n_i,&\mbox{otherwise,}
     \end{array}\right.
      &\mbox{ and }&
      
     q_i=\left\{\begin{array}{ll}
     p_i-2,&\mbox{if }p_i>p,\\
     p_i,&\mbox{if }p_i<p-1.
     \end{array}\right. 
     
 \end{array} $$

 Suppose that $\mathcal{D}_P=U_1 \cup \dots\cup U_r$, where $U_i$ is the maximum $U$-chain associated with the partition obtained in the $i$-th step of the process. Oblak conjectured that the partition $\left(|U_1|, \dots, |U_r|\right)$ is in fact the generic commuting partition $\mathcal{Q}(P)$. As we saw in Example~\ref{Uchainexa}, there may be more than one $U$-chain with maximum cardinality in each poset. Initially, Oblak addressed this issue by imposing an additional condition on the particular maximum $U$-chain that must be chosen in each step of the process. However, in \cite{K:unique}, the author proves that the restriction is not necessary and the partition described above is indeed independent of the choice of the maximum $U$-chains and is uniquely determined by $P$. Over the years since Oblak made her conjecture, results by several authors confirmed the conjecture (see \cite{Oblak}, \cite{KO}, \cite{BIK}, \cite{IK}, \cite{K:unique}, and \cite{K:smallest}). The conjecture was finally proved in its generality by Basili in \cite[Theorem 1.2]{basili21}.

\begin{figure}[!htbp]
\boxed{
\begin{tikzpicture}[scale=0.6, transform shape]
\newcommand\Square[1]{+(-#1,-#1) rectangle +(#1,#1)}
\foreach \x in {1,2,3,4,5,6}{\filldraw (\x,10-1) circle (2pt);}
\foreach \x in {1,2,3,4}{\filldraw (\x+1, 10-3) circle (2pt);}
\foreach \x in {1,2,3,4}{\filldraw (\x+1,10-4) circle (2pt);}
\foreach \x in {1,2,3}{\filldraw (\x+1.5,10-5) circle (2pt);}
\foreach \x in {1,2,3}{\filldraw (\x+1.5,10-6) circle (2pt);}
\foreach \x in {1,2}{\filldraw (\x+2,10-7) circle (2pt);}
\foreach \x in {1,2}{\filldraw (\x+2,10-8) circle (2pt);}
\foreach \x in {1}{\filldraw (\x+2.5,10-9) circle (2pt);}

\draw[black] (3.5,10-4.5) \Square{50pt};
\draw[black] (1,10-1) \Square{6pt};
\draw[black] (6,10-1) \Square{6pt};

\draw[black]
node at (3.5,0){$P=(6,\boxed{4^2,3^2},2^2,1)$}
node at (3.5,-1){$|$max U-chain$|=$ 16};

\draw[black]
node at (7.5,0){$\xrightarrow{\hspace*{1cm}}$};

\foreach \x in {2,3,4,5}{\filldraw (\x+8,10-1) circle (2pt);}
\foreach \x in {1,2}{\filldraw (\x+2+8,10-7) circle (2pt);}
\foreach \x in {1,2}{\filldraw (\x+2+8,10-8) circle (2pt);}
\foreach \x in {1}{\filldraw (\x+2.5+8,10-9) circle (2pt);}

\draw[black] (10.75,0.75) rectangle (12.25,3.25);
\draw[black] (10,9) \Square{6pt};
\draw[black] (13,9) \Square{6pt};

\draw[black]
node at (11.5,0){$P_1=(4,\boxed{2^2,1})$}
node at (11.5,-1){$|$max U-chain$|=$ 7};
\draw[black]
node at (15,0){$\xrightarrow{\hspace*{1cm}}$};

\foreach \x in {3,4}{\filldraw (\x+8+6,10-1) circle (2pt);}
\draw[black] (16.75,8.75) rectangle (18.25,9.25);
\draw[black]
node at (17.75,0){$P_2=(\boxed{2})$}
node at (17.75,-1){$|$max U-chain$|=$ 2};

\draw[black]
(0,-2-.05)--(20,-2-0.05)
(0,-2+0.05)--(20,-2+0.05);
\foreach \x in {1,2,3,4,5,6}{\filldraw (\x,10-1-12) circle (2pt);}
\foreach \x in {1,2,3,4}{\filldraw (\x+1, 10-3-12) circle (2pt);}
\foreach \x in {1,2,3,4}{\filldraw (\x+1,10-4-12) circle (2pt);}
\foreach \x in {1,2,3}{\filldraw (\x+1.5,10-5-12) circle (2pt);}
\foreach \x in {1,2,3}{\filldraw (\x+1.5,10-6-12) circle (2pt);}
\foreach \x in {1,2}{\filldraw (\x+2,10-7-12) circle (2pt);}
\foreach \x in {1,2}{\filldraw (\x+2,10-8-12) circle (2pt);}
\foreach \x in {1}{\filldraw (\x+2.5,10-9-12) circle (2pt);}

\draw[black] (4+1-3, 10-2.75-12-1-1-4+.5) rectangle (4+1, 10-2.75-12-1-1);
\draw[black] (1,10-1-12) \Square{6pt};
\draw[black] (6,10-1-12) \Square{6pt};
\draw[black] (2,10-3-12) \Square{6pt};
\draw[black] (5,10-3-12) \Square{6pt};
\draw[black] (2,10-4-12) \Square{6pt};
\draw[black] (5,10-4-12) \Square{6pt};

\draw[black]
node at (3.5,0-12){$P=(6,4^2\boxed{3^2,2^2},1)$}
node at (3.5,-1-12){$|$max U-chain$|=$ 16};
\draw[black]
node at (7.5,0-12){$\xrightarrow{\hspace*{1cm}}$};

\foreach \x in {2,3,4,5}{\filldraw (\x+8,10-1-12) circle (2pt);}
\foreach \x in {1,2}{\filldraw (\x+2+8,10-7-8) circle (2pt);}
\foreach \x in {1,2}{\filldraw (\x+2+8,10-8-8) circle (2pt);}
\foreach \x in {1}{\filldraw (\x+2.5+8,10-9-12) circle (2pt);}

\draw[black] (10.75,0.75-12) rectangle (12.25,3.25-8);
\draw[black] (10,9-12) \Square{6pt};
\draw[black] (13,9-12) \Square{6pt};

\draw[black]
node at (11.5,0-12){$P_1=(4,\boxed{2^2,1})$}
node at (11.5,-1-12){$|$max U-chain$|=$ 7};
\draw[black]
node at (15,0-12){$\xrightarrow{\hspace*{1cm}}$};

\foreach \x in {3,4}{\filldraw (\x+8+6,10-1-12) circle (2pt);}
\draw[black] (16.75,8.75-12) rectangle (18.25,9.25-12);
\draw[black]
node at (17.75,0-12){$P_2=(\boxed{2})$}
node at (17.75,-1-12){$|$max U-chain$|=$ 2};
\end{tikzpicture}
}
\caption{Performing Oblak recursive process, in two different ways, to find the generic commuting Jordan type $\mathcal{Q}(P)=(16,7,2)$ for partition $P=(6,4^2,3^2,2^2,1)$}\label{oblakfig}
\end{figure}
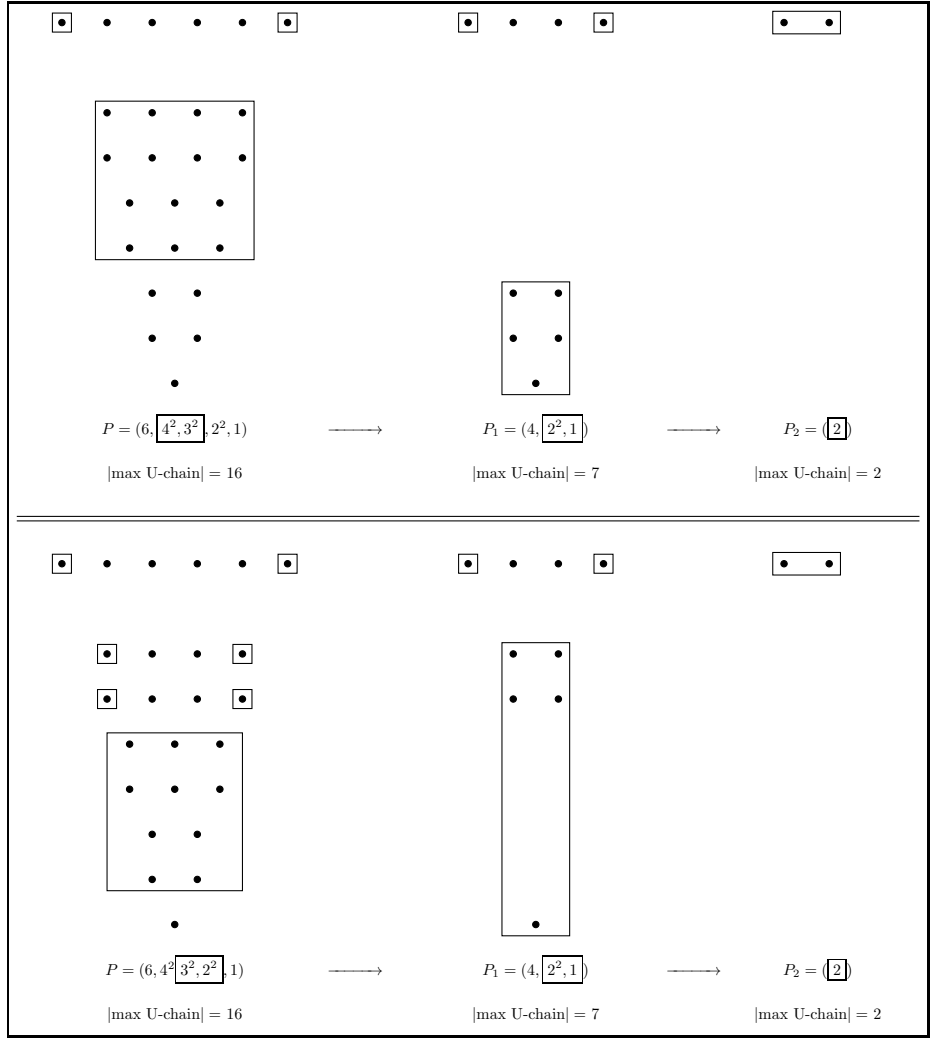

\begin{example}
    Revisiting the example $P=(6,4^2,3^3,2^2,1)$, as we saw, there are two maximum $U$-chains in $\mathcal{D}_P$, namely $U_P(4)$, corresponding to the almost rectangular subpartition $(4^2,3^2)$, and  $U_P(3)$, corresponding to the almost rectangular subpartition $(3^2,2^2)$. Both these chains have cardinality $16$.

    Assume that we start the Oblak process by removing $U_P(4)$ from $\mathcal{D}_P$. This eliminates parts of sizes 4 and 3 in $P$, shortens the part 6 of $P$ (the only one greater than 4) by 2, and keeps parts 2 and 1 (those smaller than 3) unchanged. This leads to a new partition $P_1=(4,2^2,1)$.
    Next, we perform the Oblak step for this new partition. Here, the only maximum $U$-chain is the one corresponding to the almost rectangular subpartition $(2^2,1)$. The cardinality of this $U$-chain is $2(2)+1(1)+2(1)=7$. Removing the $U$-chain from $\mathcal{D}_{P_1}$, we get partition $P_2=(2)$, in which the maximum $U$-chain (which is the whole partition) has cardinality 2. Thus the process gives rise to the partition $(16,7,2)$, which is indeed the generic partition commuting with $P$. See Figure~\ref{oblakfig}, where for simplicity, we have not included directed edges and have simply boxed vertices that correspond to the maximum $U$-chain selected in each step of the process. It is worth noting that the $U$-chain selected in each step of Oblak process is not necessarily a chain in the original poset. For example, in the example represented in Figure~\ref{oblakfig}, the maximum $U$-chain for $P_1$ is not a chain in the original poset $\mathcal{D}_P$. Recall that two elements in the poset of $P$ are comparable if and only if there is a directed path between their corresponding vertices in the directed graph $\mathcal{D}_P$ . As Figure~\ref{poset} illustrates, in the poset of $P=(6,4^2,3^2,2^2,1)$, vertices $(2,6,1)$ and $(1,2,1)$ are not comparable while they are both in the maximum $U$-chain selected in $P_1$ (the second step of the process).

    There is an alternative way to perform Oblak process for $P$, as well.  We can start the process by removing the maximum $U$-chain $U_P(3)$ from $\mathcal{D}_P$. Then continue the process as illustrated in Figure~\ref{oblakfig}. We can observe that in this example, we arrive at partition $(16,7,2)$, regardless of the choice of the maximum $U$-chain in the Oblak process.
\end{example}

\section{Partitions with the same generic commuting partition}

In previous sections, we gave an overview of results about the map $\mathcal{Q}$. In this section, we review the studies focusing of the inverse map $\mathcal{Q}^{-1}$. Recall that by Theorems~\ref{stable} and \ref{idempotent}, a partition $Q$ is in the image of the map $\mathcal{Q}$ if and only if parts of $Q$ differ pairwise by at least 2, i.e., $Q$ is stable. Fixing a stable partition $Q$, what is $\mathcal{Q}^{-1}(Q)$, the set of all partition with the same generic commuting partition $Q$? 

In \cite{Oblak12}, Oblak conjectured that if $Q=(u,u-r)$, where $r\geq 2$ and $u-r\geq 1$, then $\mathcal{Q}^{-1}(Q)$ consists of $(u-r)\times (r-1)$ partitions. She proved the conjecture for $2\leq r\leq 5$. In \cite{IKVZ}, Iarrobino, Van Steirteghem, Zhao, and the author prove the conjecture in general. Additionally, we showed that partitions in $\mathcal{Q}^{-1}(Q)$ can be arranged in a $(r-1)\times(u-r)$ table such that the partition in the $k$-th row and $\ell$-th column of the table has $k+\ell$ parts.
In \cite{IKVZ}, we also made the following general conjecture.
\begin{conjecture}(The Box Conjecture)\label{boxconj}
Let $n$ be a positive integer and $Q=(q_1, \dots, q_k)$ be a partition of $n$ such that for $q_i-q_{i+1}\geq 2$, for  $1\leq i<k$. Set 
$$s_i=\left\{\begin{array}{ll}q_i-q_{i+1}-1, &\mbox{ if }1\leq i<k,\\q_k, &\mbox{ if } i=k.\end{array}\right.$$
Then $\mathcal{Q}^{-1}(Q)$ consists of $\prod\limits_{1\leq i\leq k}s_i$ partitions. These partitions can be arranged in a $s_1\times \dots \times s_k$ ``box", $\mathcal{B}(Q)=\{P_{i_1,\dots, i_k}\,|\, 1\leq i_j\leq s_j\}$, such that $P_{i_1,\dots, i_k}$ has $\sum\limits_{1\leq j\leq k}i_j$ parts.
\end{conjecture}

In \cite{IKVZ}, the conjecture is proved for every stable partition $Q=(u+s,u,u-r)$, such that $s\leq 4$, but remains open in general.
\section{Open Questions}
We conclude this survey by stating a few related open questions for future studies.

One direction, still focusing on the generic commuting partition, is to tackle the Box Conjecture (Conjecture~\ref{boxconj}), which remains open in general. Given the nature of the statement, it is reasonable to hope for an inductive argument. In fact, partial results given in \cite{IKVZ} for the case $k=3$ solve the problem by ``decomposing" the corresponding 3-dimensional box into multiple 2-dimensional tables (sheets). However, the subtleties of the map $\mathcal{Q}(P)$, which are also captured in Oblak's process, make this approach quite involved and may require novel tools. 

Another related direction is to look beyond the generic commuting partition $\mathcal{Q}(P)$ for a given $P$ and study all Jordan types that commute with $P$. One set of partial results in this direction are given in \cite{Oblak12}, where Oblak shows the following. 

\begin{itemize}
    \item A partition $P\in \mathcal{P}(n)$ commutes with all partitions in $\mathcal{P}(n)$, if and only if $n\leq 3$, or $P=\left(2^k,1^\ell\right)$, where $k,\ell\geq 0$. 
    \item Assume that $P$ and $Q$ are both partitions of $n$ with exactly two parts. Then $P$ commutes with $Q$ if and only if $P=Q$ or $n$ is even and
$P=\left(\frac{n}{2},\frac{n}{2}\right)$ and $Q=\left(\frac{n}{2}+1,\frac{n}{2}-1\right)$.
\end{itemize}
  
  Given the fact that a pair of commuting nilpotent matrices forms an Artinian algebra, one approach to finding pairs of commuting partitions is to study Jordan types of multiplication maps in height-2 Artinian algebras. We should note that in general, the algebra generated by a pair of commuting matrices is not graded.
\bibliography{ref}
\bibliographystyle{amsplain}
\end{document}